\documentclass[11pt]{article}
\usepackage{amsmath}
\usepackage{amssymb}
\usepackage{amscd}
\usepackage[all]{xy}
\usepackage[T1]{fontenc}
\usepackage{epsfig}
\newtheorem{theo}{Theorem}[section]
\newtheorem{prop}[theo]{Proposition}

\newtheorem{lemm}[theo]{Lemma}

\title{\textbf{\textsc{free quantum analogues of the 
first fundamental theorems of invariant theory}}} 
\author{\textsc{Julien Bichon}}
\date{{\small \textsl{Laboratoire de Math\'ematiques Appliqu\'ees,
Universit\'e de Pau et des Pays de l'Adour, \\
IPRA, Avenue de l'universit\'e,
64000 Pau, France.}
E-mail: Julien.Bichon@univ-pau.fr}}

\makeatletter
\renewcommand{\@makefnmark}{}
\makeatother

\begin{document}

\maketitle

\begin{abstract}
We formulate and prove a free quantum analogue of the
first fundamental theorems of invariant theory.
More precisely, the polynomial function algebras
on matrices are replaced by free algebras, while
the universal cosovereign Hopf algebras play the role of 
the general linear group. 
\end{abstract}

\section*{Introduction}

This note was inspired by reading the paper \cite{[GLR]}
of Goodearl, Lenagan and Rigal.  
We have followed part of their introduction very closely and
we have tried to keep, as much as possible, the terminology
and notation of this paper. 

\medskip

Consider a fixed field $K$ and positive integers $m,n,t>0$.
For integers $u,v>0$, we write $M_{u,v}$ for the set of
$u \times v$ matrices with entries in $K$. The general
linear group $GL_t = GL_t(K)$ acts on the variety 
$V = M_{m,t} \times M_{t,n}$, with action given by
\begin{align*}
GL_t & \times V 
\longrightarrow V \\
(g, (A,& B)) \longmapsto (Ag^{-1}, gB).
\end{align*} 
Thus $GL_t$ acts on the algebra 
$\mathcal O(V) \cong \mathcal O(M_{m,t}) \otimes \mathcal (M_{t,n})$.
The first fundamental theorems of invariant theory
describe the subalgebra 
$\mathcal O(V)^{GL_t}$ of invariants for this action in the following
way. Consider the multiplication map
\begin{align*}
\theta : M_{m,t} \times & M_{t,n} \longrightarrow M_{m,n}\\
(A, & B) \longmapsto AB
\end{align*}
and denote by $\theta^* : \mathcal O(M_{m,n}) \longrightarrow
\mathcal O(M_{m,t}) \otimes \mathcal O(M_{t,n})$ the induced algebra
morphism. 

\begin{theo}
The ring of invariants $\mathcal O(V)^{GL_t}$ equals 
{\rm Im}$\theta^*$.
\end{theo} 

\begin{theo}
Let $X_{ij}$ ($1 \leq i \leq m$ , $1 \leq j \leq n$)
be the usual coordinate functions on $M_{m,n}$. Let 
$\mathcal T_{t+1}$ be the ideal of $\mathcal O(M_{m,n})$
generated by the $(t+1) \times (t+1)$ minors of the matrix
$(X_{ij})$ over $\mathcal O(M_{m,n})$ (this ideal is zero if
$t \geq {\rm min} \{m,n\}$). Then the kernel of 
$\theta^*$ is $\mathcal T_{t+1}$.
\end{theo}

These two theorems are respectively known as the first and second fundamental
theorems of invariant theory (for $GL_t$). See \cite{[DCP]}.
They give a full
description of the algebra $\mathcal O(V)^{GL_t}$.

\medskip

Goodearl, Lenagan and Rigal \cite{[GLR]} generalized these theorems
for quantized coordinate algebras. In this paper we prove analogues
of these theorems for free algebras.

Let $u,v >0$ be positive integers, and denote by $A(u,v)$
the free algebra on $uv$ generators. The natural analogue of the 
algebra morphism $\theta^*$ above is the algebra morphism
(denoted by $\theta$ for simplicity)
\begin{align*}
\theta :  A(m & ,n) \longrightarrow A(m,t) \otimes A(t,m)\\
& x_{ij} \longmapsto \sum_{k=1 }^t y_{ik} \otimes z_{kj}
\end{align*} 
where $x_ {ij}$, $y_{ij}$, $z_{ij}$ stand for generators of
$A(m,n)$, $A(m,t)$ and $A(t,n)$ respectively.
It is easy to show that contrary to the commutative case,
the morphism $\theta$ is always injective. The following question arises
naturally: does there exist a quantum group $G$ acting on $A_{m,t}$ 
and $A_{t,n}$ and such that Im$\theta$ equals 
$(A_{m,t} \otimes A_{t,n})^G$ ? 
We show that the universal cosovereign Hopf algebras introduced
in $\cite{[Biso]}$, which are natural
free analogues of the general linear groups in quantum group theory,
 answer positively to this question.
The key for proving this result is a representation theoretic
property of these Hopf algebras, proved in \cite{[Biun]}.

\medskip

Our work is organized as follows. In section 1 we recall the 
setup of \cite{[GLR]}
for non-commutative analogs of the first fundamental 
theorems of invariant theory.
In Section 2 we recall some basic facts concerning the universal
cosovereign Hopf algebras and state the main theorem. The proof
is given in Section 3.

\medskip

Throughout the paper we work over an arbitrary base field $K$.

\section{The setup for non-commutative invariant theory}

Let us first recall the setup, due to Goodearl,
Lenagan and Rigal, 
for stating non-commutative analogues of the first fundamental
theorems of invariant theory. Similar considerations were done
independently by Banica \cite{[Ba]} in the context of Kac algebras
actions.

\medskip

Let $H$ be a Hopf algebra, let $(A, \rho)$ be a right $H$-comodule algebra
and let $(B, \lambda)$ be a left $H$-comodule algebra, where
$$\rho : A \longrightarrow A \otimes H \quad {\rm and} \quad
\lambda : B \longrightarrow H \otimes B$$
are the coactions of $H$.
One can turn $A$ into a left $H$-comodule with the coaction
$\rho' = \tau \circ ({\rm Id}_A \otimes S) \circ \rho$
where $\tau : A \otimes H \longrightarrow H \otimes A$ is the standard flip.
Thus one can consider the tensor product $A \otimes B$
of the left $H$-comodules $A$
and $B$ which is a left $H$-comodule, but which, 
in the non-commutative situation,
is not any more a  comodule algebra in general.
However, we have the following result from \cite{[GLR]}
(see also \cite{[Ba]}):

\begin{prop}
The set of coinvariants $(A \otimes B)^{{\rm co}H}$ is a subalgebra
of $A \otimes B$.
\end{prop}

Here it should be recalled that if $M$ is a left $H$-comodule
with coaction $\alpha : M  \longrightarrow H \otimes M$, the set of 
coinvariants is 
$M^{\rm co H} = \{x \in M \ | \ \alpha(x) = 1 \otimes x \}$. 
Denoting by $\mathcal C$ the category of left $H$-comodules
and by $I$ the trivial one-dimensional comodule, then
$M^{\rm co H}$ is canonically identified with 
Hom$_{\mathcal C}(I, M)$.

\medskip

Fix positive integers $m,n,t>0$. 
We denote by $A(m,n)$, $A(m,t)$ and $A(t,n)$ the free algebras on 
$mn$, $mt$ and $tn$ generators respectively, with
canonical generators denoted respectively by $x_{ij}$, $y_{ij}$ and $z_{ij}$. 
The free analogue of the comultiplication map $\theta^*$ of the introduction
is the  algebra morphism
\begin{align*}
\theta :  A(m & ,n) \longrightarrow A(m,t) \otimes A(t,m)\\
& x_{ij} \longmapsto \sum_{k=1 }^t y_{ik} \otimes z_{kj}.
\end{align*} 
It may be shown easily by direct computations that
$\theta$ is injective. 

\smallskip

Now consider a Hopf algebra $H$ having a multiplicative matrix 
$u = (u_{ij}) \in M_t(H)$. 
This means that for $i,j \in \{1, \ldots , t \}$, we have
$\Delta(u_{ij}) = \sum_k u_{ik} \otimes u_{kj}$ and 
$\varepsilon(u_{ij}) = \delta_{ij}$.
Then $A(m,t)$ is a right
$H$-comodule algebra with coaction
\begin{align*}
\rho : A(m & ,t) \longrightarrow A(m,t) \otimes H\\
& y_{ij} \longmapsto \sum_{k=1}^t y_{ik} \otimes u_{kj}.
\end{align*} 
In the same way $A(t,n)$ is a left $H$-comodule algebra
with coaction
\begin{align*}
\lambda : A(t & ,n) \longrightarrow H \otimes A(t,n)\\
& z_{ij} \longmapsto \sum_{k=1}^t u_{ik} \otimes z_{kj}.
\end{align*} 
Thus we are in the preceding situation and we can consider
the algebra $(A(m,t) \otimes A(t,n))^{{\rm co} H}$.
Similarly to Proposition 2.3 in \cite{[GLR]}, we have
the following result.

\begin{prop}
The algebra $(A(m,t) \otimes A(t,n))^{{\rm co} H}$
is a subalgebra of $A(m,t) \otimes A(t,n)$ containing
${\rm Im}\theta$.
\end{prop}

It is natural to wonder when the algebra morphism 
$\theta$ is surjective. We will see that $\theta$ is surjective
for the universal cosovereign Hopf algebras \cite{[Biso]}.
In fact the key property in order that $\theta$ be surjective
is that the tensor powers of the comodule $U$ associated to the 
multiplicative matrix $u$ be simple non-equivalent comodules.
At the ring-theoretic level, this corresponds to the fact that
the elements $u_{ij}$ generate a free subalgebra on $t^2$ generators.

\section{Universal cosovereign Hopf algebras and the theorem} 

Let $F \in GL_t$.
Recall \cite{[Biso]} that 
the algebra $H(F)$ is defined to be the universal algebra
with generators
 $(u_{ij})_{1 \leq i,j \leq t}$,
 $(v_{ij})_{1 \leq i,j \leq t}$ and relations:
$$ {u} {^t \! v} = {^t \! v} u = I_t = {vF} {^t \! u} F^{-1} = 
{F} {^t \! u} F^{-1}v \ ,$$
where $u= (u_{ij})$, $v = (v_{ij})$ and $I_t$ is
the identity $t \times t$ matrix. It turns out \cite{[Biso]} that
$H(F)$ is a Hopf algebra
with comultiplication defined by 
$\Delta(u_{ij}) = \sum_k u_{ik} \otimes u_{kj}$
and $\Delta(v_{ij}) = \sum_k v_{ik} \otimes v_{kj}$, with counit
defined by $\varepsilon (u_{ij}) = \varepsilon (v_{ij}) = \delta_{ij}$ and
with antipode defined by $S(u) = {^t \! v}$ and $S(v) = F {^t \! u} F^{-1}$.
Furthermore $H(F)$ is a cosovereign Hopf algebra \cite{[Biso]}: there
exists an algebra morphism $\Phi : H(F) \longrightarrow k$
such that $S^2 = \Phi * {\rm id} * \Phi^{-1}$.
 The Hopf algebras 
$H(F)$ have the following universal property (\cite{[Biso]}, Theorem 3.2).

\smallskip

\textsl{Let $H$ be Hopf algebra and let $V$ be finite dimensional
$H$-comodule isomorphic with its bidual comodule $V^{**}$. 
Then there exists a matrix $F \in GL_t$ ($t = \dim V$) such that
$V$ is an $H(F)$-comodule and such that there exists a Hopf algebra
morphism $\pi : H(F) \longrightarrow H$ such that 
$(1_V \otimes \pi) \circ \beta_V = \alpha_V$,
where $\alpha_V : V \longrightarrow V \otimes H$ and
$\beta_V : V \longrightarrow V \otimes H(F)$ denote the coactions 
of $A$ and $H(F)$ on $V$ respectively.
In particular 
every finite type cosovereign Hopf algebra is a homomorphic quotient
of a Hopf algebra $H(F)$.}

\smallskip

In view of this universal property it is natural to say that the 
Hopf algebras $H(F)$ are the universal cosovereign Hopf algebras,
or the free cosovereign Hopf algebras, and to see these Hopf algebras
as natural analogues of the general linear groups in quantum
group theory.

\medskip

Coming back to the situation of the first section,
we have the following result, which is a free quantum analogue
of the first fundamental theorems of invariant theory.   

\begin{theo}
Let $m,n,t>0$ be positive integers and 
let $F \in GL_t$. Then the algebra morphism
$$\theta : A(m,n) \longrightarrow (A(m,t) \otimes A(t,n))^{{\rm co} H(F)}$$
is an isomorphism.
\end{theo}

\section{Proof of the theorem}

\noindent
\textbf{3.1.}
We begin with some general considerations. 
Let $\mathcal C$ be a $K$-linear strict
tensor category, that is $\mathcal C$ is an abelian $K$-linear category
and $\mathcal C = (\mathcal C , \otimes , I)$
is a strict tensor category \cite{[Ka]} such the 
tensor product is $K$-linear in each variable.
Let $X,Y$ be some objects of $\mathcal C$. We define a $K$-algebra
$\mathcal C(X,Y)$ in the following way. As a vector space
$$\mathcal C(X,Y) = \bigoplus_{k \in \mathbb N} {\rm Hom}_{\mathcal C}
(X^{\otimes k} , Y^{\otimes k}).$$
The product of $\mathcal C(X,Y)$ is defined on homogeneous elements
$f_1 \in {\rm Hom}_{\mathcal C} (X^{\otimes k_1} , Y^{\otimes k_1})$ and
$f_2 \in {\rm Hom}_{\mathcal C} (X^{\otimes k_2} , Y^{\otimes k_2})$
by $f_1f_2 := f_1 \otimes f_2 \in 
{\rm Hom}_{\mathcal C} (X^{\otimes k_1 +k_2} , Y^{\otimes k_1+k_2})$.
It is immediate to check that $\mathcal C(X,Y)$ is an
associative $K$-algebra.

Let $m,n \in \mathbb N^*$ and let $U$ be an object of $\mathcal C$. 
Let $X= U^m$ be the direct sum of $m$ copies of $U$. Let 
$p_i : U^m \longrightarrow U$ and $v_i : U \longrightarrow U^m$,
$1 \leq i \leq m$, the canonical morphisms such that
$\sum_{i=1}^m v_i \circ p_i = 1_X$ and $p_i \circ v_j = \delta_{ij}1_U$.
Similarly, let $Y= U^n$ be the direct sum of $n$ copies of $U$. Let 
$q_i : U^n \longrightarrow U$ and $u_i : U \longrightarrow U^n$,
$1 \leq i \leq n$, the canonical morphisms such that
$\sum_{i=1}^m u_i \circ q_i = 1_Y$ and $q_i \circ u_j = \delta_{ij}1_U$.

\begin{lemm}
Assume that ${\rm End}_{\mathcal C}(U) = K$. Then the algebra
morphism
\begin{align*}
\psi : A(m & ,n) \longrightarrow \mathcal C(X,Y) \\
& x_{ij} \longmapsto u_j \circ p_i
\end{align*}
is an isomorphism
\end{lemm}

\noindent
\textbf{Proof}. Let $i_1, \ldots , i_k \in \{1 , \ldots , m\}$,
$j_1, \ldots , j_k \in \{1 , \ldots , m\}$. Then
$$\psi(x_{i_1j_1} \ldots x_{i_k j_k}) =
(u_{j_1} \otimes \ldots \otimes u_{j_k}) \circ
(p_{i_1} \otimes \ldots \otimes p_{i_k}).$$
Thus $\psi$ transforms a basis of $A(m,n)$ into a basis
of $\mathcal C(X,Y)$, and is an isomorphism. $\square$

\bigskip

\noindent
\textbf{3.2.} Consider again the $K$-linear strict tensor category
$\mathcal C$. Let $X$ be an object of $\mathcal C$. Recall
\cite{[Ka]} that a right dual for $X$ is a triplet
$(X^*,e,d)$ where $X^*$ is an object of $\mathcal C$, while
$e : X \otimes X^* \rightarrow I$ and $d : I \rightarrow X^* \otimes X$
are morphisms such that
$$ 
(e \otimes 1_X) \circ (1_X \otimes d) = 1_X \quad {\rm and}
\quad (1_{X^*} \otimes e) \circ (d \otimes 1_{X^*}) = 1_{X^*}. 
$$
We then have for all objects $Y,Z$ of $\mathcal C$ isomorphisms
\begin{align*}
{\rm Hom}_{\mathcal C} (Y, X^* & \otimes Z) \cong
{\rm Hom}_{\mathcal C}(X \otimes Y, Z) \\
& f \ \longmapsto \ (e \otimes 1_Z) \circ (1_X \otimes f).
\end{align*}
Also recall that if $X$ has a right dual $(X^*,e,d)$, then 
$X^{\otimes n}$ has a right dual $(X^{*\otimes n},e_n,d_n)$
for all
$n \in \mathbb N^*$, where
$e_n= e \circ \ldots \circ (1_{X^{\otimes n-1}} \otimes e \otimes
1_{X^{*\otimes n-1}})$ and
$d_n = ( 1_{X^{*\otimes n-1}} \otimes d \otimes 1_{X^{\otimes n-1}})
\circ \ldots \circ d$.

Recall finally that in the tensor category of finite-dimensional
left comodules over a Hopf algebra, every object has a right
dual.

\bigskip

\noindent
\textbf{3.3.} Let $m,n,t >0$ be positive integers and let $F \in GL_t$.
Consider the multiplicative matrix $u =(u_{ij}) \in M_t(H(F))$.
One associates a right $H(F)$-comodule $U_r$ to $u$:
as a vector space $U_r = K^t$ with its standard basis 
$e_1, \ldots , e_t$, 
 and the coaction $\alpha : U_r \longrightarrow U_r \otimes H(F)$
is given by $\alpha(e_i) = \sum_j e_j \otimes u_{ji}$.
Similarly one associates a left $H(F)$-comodule $U_l$ to $u$ :
as a vector space $U_l = K^t$ and the coaction
$\beta : U_l \longrightarrow H(F) \otimes U_l$ is given by
$\beta(e_i) = \sum_j u_{ij} \otimes e_j$.

We have the following key result from \cite{[Biun]} (Corollary 2.6):

\begin{prop}
The left (respectively right) $H(F)$-comodules
$U_l^{\otimes k}$, $k \in \mathbb N^*$ (respectively
$U_r^{\otimes k}$, $k \in \mathbb N^*$) are simple
non-equivalent left (respectively right) $H(F)$-comodules, and have
endomorphism algebras isomorphic with $K$.
\end{prop}

As a right $H(F)$-comodule algebra $A(m,t)$ is naturally identified
with the tensor algebra
$$T(U_r^m) = \bigoplus_{i \in \mathbb N} (U_r^m)^{\otimes i}.$$ 
Now transform $A(m,t)$ and $T(U_r^m)$ into left $H(F)$-comodules
in the manner of Section 1. Then it is immediate that, as a left
$H(F)$-comodule, $T(U_r^{m})$ is identified with 
$T(U_l^{m*})$. We have an isomorphism
\begin{align*}
\phi_1 :  A(m & ,t) \longrightarrow T(U_l^{m*})^{\rm op} \\
& y_{ij} \longmapsto v_i(e_j)^*
\end{align*}
which is both an algebra isomorphism and an $H(F)$-comodule isomorphism
(we use the notations of 3.1).
Similarly $A(t,n)$ is identified with $T(U_l^n)$ by the following
left $H(F)$-comodule algebra isomorphism:
\begin{align*}
\phi_2 : A(t & ,n) \longrightarrow T(U_l^{n}) \\
& z_{ij} \longmapsto u_j(e_i).
\end{align*}

We have now all the ingredients to prove Theorem 2.1.
We work in the $K$-linear tensor category $\mathcal C$ of
left $H(F)$-comodules. Since $\mathcal C$ is a concrete tensor 
category of vector spaces, we can proceed as if it was strict and
the considerations of 3.1 and 3.2 are valid. 
We put $U = U_l$. We have

\begin{align*}
( & T(U^{m*}) \otimes T(U^n))^{{\rm co}H(F)} 
\cong {\rm Hom}_{\mathcal C}(I, T(U^{m*}) \otimes T(U^n)) \\
& \cong {\rm Hom}_{\mathcal C}(I, \bigoplus_{i,j \in \mathbb N}
(U^{m*})^{\otimes i} \otimes (U^n)^{\otimes j})
\cong  \bigoplus_{i,j \in \mathbb N} 
{\rm Hom}_{\mathcal C}(I, (U^{m*})^{\otimes i} \otimes (U^n)^{\otimes j})\\
& \cong  \bigoplus_{i,j \in \mathbb N} 
{\rm Hom}_{\mathcal C}((U^{m})^{\otimes i} , (U^n)^{\otimes j})
\quad ({\rm by} \ 3.2) \\
& \cong \bigoplus_{i \in \mathbb N} 
{\rm Hom}_{\mathcal C}((U^{m})^{\otimes i},(U^n)^{\otimes i}) \quad
({\rm by} \ {\rm Proposition} \ 3.2) \ = \mathcal C(U^m,U^n) .
\end{align*}

It is straightforward to check that for
$i_1, \ldots , i_p \in \{1, \ldots , m \}$ and 
$j_1, \ldots , j_p \in \{1, \ldots , n \}$, the
isomorphism just considered transforms the element
$$
\sum_{k_1, \ldots , k_p}
v_{i_p}(e_{k_p})^* \otimes \ldots \otimes v_{i_1}(e_{k_1})^*
\otimes u_{j_1}(e_{k_1}) \otimes \ldots \otimes u_{j_p}(e_{k_p})$$
into the element
$$
(u_{j_1} \otimes \ldots \otimes u_{j_k}) \circ
(p_{i_1} \otimes \ldots \otimes p_{i_k})$$
of $\mathcal C(U^m,U^n)$. Hence this isomorphism, composed
with $(\phi_1 \otimes \phi_2) \circ \theta$, yields the 
isomorphism $\psi$ of Lemma 3.1, and
$\theta$ is itself an isomorphism. This concludes the proof of Theorem 2.1.

\bigskip

\noindent
\textbf{Remark}. 
Theorem 2.1 is still valid if one replaces the Hopf algebra
$H(F)$ by 
$H(t)$, the free Hopf algebra generated by the matrix coalgebra
$M_t^*$ \cite{[T]}. 
This is clear from the proof and from the fact
that $H(F)$ being a quotient of 
$H(t)$, Proposition 3.2 remains valid for $H(t)$.


\begin{thebibliography}{25}

\small{

\bibitem{[Ba]} T. \textsc{Banica}, Subfactors associated
to compact Kac algebras, Integral Equations Oper. Theory 39
(2000), No.1, 1-14.

\bibitem{[Biso]} J. \textsc{Bichon}, Cosovereign Hopf algebras,
J. Pure Appl. Algebra 157, No.2-3 (2001), 121-133.

\bibitem{[Biun]} J. \textsc{Bichon},
Corepresentation theory of universal cosovereign Hopf algebras,
Preprint, 2002.

\bibitem{[DCP]} C. \textsc{De Concini}, C. \textsc{Procesi},
A characteristic free approach to invariant theory,
Adv. Math. 21 (1976), 330-354.

\bibitem{[GLR]} K.R. \textsc{Goodearl}, T.H. \textsc{Lenagan},
L. \textsc{Rigal}, The first fundamental theorem of 
coinvariant theory for the quantum general linear group, 
Publ. Res. Inst. Math. Sci. 36 (2000), No.2, 269-296.

\bibitem{[Ka]} C. \textsc{Kassel}, \textsl{Quantum groups}, 
GTM 155, Springer, 1995

\bibitem{[T]} M. \textsc{Takeuchi}, Free Hopf algebras generated
by coalgebras, J. Math. Soc. Japan 23, No.4 (1971), 581-562.

}

\end{thebibliography}
\end{document}